\documentclass[10pt]{amsart}
\usepackage{amssymb,latexsym}
\theoremstyle{plain}
\newtheorem{theorem}{Theorem}

\newtheorem*{TA}{Theorem A}
\newtheorem*{TB}{Theorem B}
\newtheorem*{TC}{Theorem C}

\theoremstyle{definition}

\theoremstyle{remark}

\numberwithin{equation}{section}

\newcommand{\R}{\mathbb R}

\begin{document}

\title[The Benjamin-Ono equation]{The IVP for the Benjamin-Ono equation in weighted Sobolev spaces II}
\author{Germ\'an Fonseca}
\address[G. Fonseca]{Departamento  de Matem\'aticas\\
Universidad Nacional de Colombia\\ Bogota\\Colombia}
\email{gefonsecab@unal.edu.co}
\author{Felipe Linares}
\address[F. Linares]{Instituto de Matematica Pura e Aplicada\\
IMPA\\Rio de Janeiro\\Brazil}
\email{linares@impa.br}

\author{Gustavo Ponce}
\address[G. Ponce]{Department  of Mathematics\\
University of California\\
Santa Barbara, CA 93106\\
USA.}
\email{ponce@math.ucsb.edu}
\begin{abstract} 
In this work we continue our study initiated in \cite{GFGP} on the uniqueness properties 
of real solutions to the IVP associated to the Benjamin-Ono (BO) equation.
In particular, we shall show that the uniqueness results established in \cite{GFGP} do not extend 
to any  pair of non-vanishing solutions of the BO equation. Also, we shall prove that the uniqueness
result established in \cite{GFGP} under a hypothesis involving information of the solution at three 
different times can not be relaxed to  two different times.
\end{abstract}
\maketitle
\begin{section}{Introduction}\label{S: 1}

This work is concerned with special decay and  uniqueness properties 
of real solutions of the initial value problem (IVP)
for the  Benjamin-Ono (BO) equation  
\begin{equation}\label{BO}
\begin{cases}
\partial_t u + \mathcal  H\partial_x^2u +u\partial_x u = 0, \qquad t, x
\in \R,\\
u(x,0) = u_0(x),
\end{cases}
\end{equation}
where $ \mathcal  H$ denotes  the Hilbert transform
\begin{equation}
\label{hilbertt}
\begin{aligned}
\mathcal H f(x)&=\frac{1}{\pi} \,\text{p.v.} (\frac{1}{x}\ast f)(x)\\
&=\frac{1}{\pi}\lim_{\epsilon\downarrow 0}\int_{|y|\geq \epsilon}\frac{f(x-y)}{y}dy=-i\,(\text{sgn}(\xi)\,\widehat {f}(\xi))^{\lor}(x).
\end{aligned}
\end{equation}

The BO equation  was deduced by Benjamin
\cite{Be} and Ono \cite{On} as  a model for long internal gravity waves in deep stratified fluids.
Later, it was also shown that it is  a completely integrable system (see \cite{AbFo}, \cite{CoWi} and references therein).
 
The problem of finding the minimal regularity property (measured in the Sobolev scale 
$H^s(\R)= \left(1-\partial^2_x\right)^{-s/2} L^2(\R),\,s\in\R$) required
in the data $u_0$ which guarantees that the IVP \eqref{BO} is locally wellposed (LWP) 
(global wellposed (GWP) follows by combining the LWP and the conservation laws) 
has been extensively considered. Thus, one has the following list of works:
in \cite{Sa} $s>3$ was proven, in \cite{ABFS} and 
\cite{Io1} $s>3/2$, 
in \cite{Po} $s\geq 3/2$, 
in \cite{KoTz1}  $s>5/4$, 
in \cite{KeKo} $s>9/8$,
in \cite{Ta} $s\geq 1$,
in \cite{BuPl} $s>1/4$,
and 
in \cite{IoKe} $s\geq 0$. It should be pointed out that the result in \cite{MoSaTz} (see also \cite{KoTz2}) 
implies that no well-posedness result in $H^s(\R),\;s\in\R\,$ can be established by a solely contraction principle arguments.
For further results on uniqueness and comments we refer to \cite{LMDP}. 

Our study here includes both the regularity and the decay of the solution measured in the $L^2$ sense. 
More precisely, we deal with persistence properties (i.e. if the data $u_0$ belongs to the function space $X$, 
then the corresponding solution 
of \eqref{BO} defines a continuous curve on $X$, $u\in C([0,T]:X)$) of real valued solutions of the IVP \eqref{BO}  
in the weighted Sobolev spaces
\begin{equation}
\label{spaceZ}
Z_{s,r}=H^s(\R)\cap L^2(|x|^{2r}dx),\;\;\;\;\;\;s,\,r\in\R,
\end{equation}
and 
\begin{equation}
\label{spaceZdot}
\dot Z_{s,r}=\{ f\in H^s(\R)\cap L^2(|x|^{2r}dx)\,:\,\widehat {f}(0)=0\},\;\;\;\;\;\;\;s,\,r\in\R.
\end{equation}
 Notice that the conservation law for  solutions of \eqref{BO} 
 $$
 I_1(u_0)=\int_{-\infty}^{\infty} u_0(x)dx=\int_{-\infty}^{\infty} u(x,t)dx,
 $$
 guarantees that the property $\,\widehat{u}_0(0)=0$
 is preserved by the solution flow.
 
  As an extension of the results in \cite{Io1}, \cite{Io2}  the following  
theorems were proven in \cite{GFGP}:
    
 \begin{TA}   $($\cite{GFGP}$)$ \label{theorem3} 
  (i) Let $s\geq 1, \;r\in [0,s]$, and $\,r<5/2$. If $u_0\in Z_{s,r}$, then the solution $u(x,t)$ of the IVP 
\eqref{BO} satisfies that $ u\in C(\R:Z_{s,r})$.
 
 (ii) For  $s>9/8$  ($s\geq 3/2$), $\;r\in [0,s]$, and $\,r<5/2$ the IVP \eqref{BO} is LWP (GWP resp.) in $Z_{s,r}$.
 
(iii) If  $\,r\in [5/2,7/2)$ and $\,r\leq s$, then the IVP \eqref{BO} is GWP 
in $\dot Z_{s,r}$.

\end{TA}

\begin{TB}   $($\cite{GFGP}$)$ \label{theorem4} 
 Let $u\in C(\R : Z_{2,2})$ be a solution of the IVP \eqref{BO}. If   there exist  two different times
 $\,t_1, t_2\in \R$ such that 
 \begin{equation}
 \label{2timesw} 
 u(\cdot,t_j)\in Z_{5/2,5/2},\;\;j=1,2,\;\;\text{then}\;\;\;\widehat {u}_0(0)=0\,,\;\;
(\text{so}\;\; u(\cdot, t)\in  \dot Z_{5/2,5/2},\,\, \forall t\in R).
 \end{equation}
 
\end{TB}
 
 \vskip.1in
 
 \begin{TC}   $($\cite{GFGP}$)$ \label{theorem5} 
 Let $u\in C(\R : \dot Z_{3,3})$ be a solution of the IVP \eqref{BO}. If   there exist  three different times
 $\,t_1, t_2, t_3\in \R$ such that 
 \begin{equation}
 \label{3times} 
 u(\cdot,t_j)\in Z_{7/2,7/2},\;\;j=1,2,3,\;\;\text{then}\;\;\;u(x,t)\equiv 0.
 \end{equation}
 
\end{TC}
\vskip.1in

\underline{Remarks} : (a) Theorem  A part (ii) with $s\geq r=2$, Theorem A part (iii) with $s\geq r=3$, and Theorem C with $s=r=4$ were 
proved by Iorio, see \cite{Io1}, \cite{Io2}.

(b) Theorem B  shows that the condition $\widehat {u}_0(0)=0$ is necessary for  persistence property
of the solution to hold in  $Z_{s,5/2}$,  with $s\geq 5/2$, so in that sense 
of Theorem A parts (i)-(ii) are optimal. 
Theorem C affirms  that there is 
an upper limit of the spatial $L^2$-decay rate of the solution i.e.  
$$
|x|^{7/2}u(\cdot,t)\notin L^{\infty}([0,T]: L^2(\R)), \;\;\;\;  \,\,\,\,\forall \,T>0,
$$
regardless of the decay and  regularity of the non-zero initial data $u_0$. 
In particular, Theorem C shows that Theorem A part (iii) is sharp.

\vskip.1in

In view of the results in Theorem A, Theorem B, and Theorem C the following two questions present themselves.
\vskip.1in
\underline{Question 1 } : Can these uniqueness results  be extended to any pair of solutions $\,u_1,\;u_2\,$
of the \eqref{BO} with $u_1\neq 0,\;u_2\neq 0$?
\vskip.1in
We recall that the uniqueness results obtained in \cite{EKPV1} for the  IVP associated to the $k$-generalized Korteweg-de Vries ($k$-gKdV) equation
\begin{equation}\label{kgKdV}
\partial_t u +\partial_x^3u +u^k\partial_x u = 0, \qquad t, x
\in \R,
\end{equation}
and those in \cite{EKPV2} for the IVP associated to the semi-linear Schr\"odinger (NLS) equation 
\begin{equation}
\label{NLS}
i\partial_t u +\Delta u = F(u,\bar u), \qquad t\in \R, x
\in \R^n,
\end{equation}
hold for any pair $u_1,\,u_2\,$ of solutions in a suitable class.

  Our first result gives a negative answer to Question 1:

\begin{theorem}\label{theorem6} 
 There exist $\,u_1,\;u_2\in C(\R : Z_{4,2}),\;u_1\neq 0,\,u_2\neq 0$  solutions of the IVP \eqref{BO} such that  
  $$
 u_1\neq u_2
 $$
 and for any $T>0$
 \begin{equation}
 \label{2atimeus} 
 u_1-u_2\in L^{\infty}([-T,T] : Z_{4,4}).
 \end{equation}
 \end{theorem}
\vskip.1in

\underline{Remarks} : (a)  Combining the argument presented here and those used in \cite{GFGP} 
relying on the notion of $A_2$ weight one can extend the result in Theorem 
\ref{theorem6} to the index $Z_{5,9/2-}$ by assuming that $\,u_1,\;u_2\in C(\R : Z_{5,2})$. 
This tells us that for  a uniqueness result compromising  any pair of suitable solutions of the BO equation
to be valid it should require a decay index $r\geq 9/2$.

\vskip.1in

Next, we observe that  the hypothesis \eqref{3times} in  Theorem C 
involves a condition  on the solution $u(x,t)$ at three different times $\,t_1<t_2<t_3$.

\vskip.1in
\underline {Question 2}:  Can the assumption  \eqref{3times} in Theorem C be reduced to  two different times $\,t_1<t_2$?
\vskip.1in
We recall that the uniqueness results for the $k$-gKdV  in \cite{EKPV1}, for the NLS in \cite{EKPV2}, those obtained in
\cite{HMPZ} for the IVP associated to the Camassa-Holm equation 
\begin{equation}
\label{CH}
\partial_t u-\partial_{t}\partial^2_{x}u+3u\partial_x u - 2\partial_x u \partial^2_x u - u\partial_x^3 u= 0,\quad t, x
\in \R,
\end{equation}
as well as  many other deduced for dispersive models  require a condition involving only two different times. 
\vskip.1in

Surprisingly, our second result shows that for the BO this is not the case, the condition involving three different 
times in Theorem C is necessary:
\vskip.1in
$\exists \;u\in C(\R: \dot Z_{5,7/2-}),\;u\ne 0,$ solution of \eqref{BO} for which there are $t_1,\,t_2\in\R,\;t_1\neq t_2$ 
such that
$$
u(\cdot,t_j)\in \dot Z_{5,4}\subset \dot Z_{7/2,7/2-},\;\;\;\;\;\;j=1,2.
$$
\vskip.1in
More precisely, we shall prove :

\begin{theorem}
\label{theorem7} 
 For any $u_0\in \dot Z_{5,4}$ such that
 \begin{equation}
 \label{cond1}
  \int_{-\infty}^{\infty} x \,u_0(x)dx\neq 0,
 \end{equation}
the corresponding solution $\,u\in C(\R : \dot Z_{5,7/2-})$  of the \eqref{BO} provided by Theorem A part (iii) 
satisfies that
 \begin{equation}
 \label{attime} 
 u(\cdot,t^*)\in \dot Z_{5,4},
 \end{equation}
 where
  \begin{equation}
 \label{thetime} 
 t^*= - \frac{4}{\|u_0\|_2^2}\;\int_{-\infty}^{\infty} x\, u_0(x)dx.
 \end{equation}
 \end{theorem} 
 \vskip.1in

\underline{Remarks} : (a) The result in Theorem \ref{theorem7} is due to the relationship between the dispersive 
part  and the structure of the nonlinearity
of the BO equation. In particular, one can see that if $u_0\in \dot Z_{5,4}$ verifies  \eqref{cond1},
 then the solution $\,W(t)u_0(x)\,$of the associated linear IVP 
$$
\partial_t u + \mathcal  H\partial_x^2u=0,\;\;\;\;\;u(x,0)=u_0(x),
$$
satisfies 
$$
W(t)u_0(x)=c(e^{-it|\xi|\xi}\widehat {u_0}(\xi))^{\lor}\in L^2(|x|^{7-})- L^2(|x|^7),\;\;\;\;\;\;\forall \,t\neq 0.
$$
However, for the same data $u_0$ one has that the solution $u(x,t)$ of the IVP \eqref{BO} satisfies
$$
u(\cdot,0),\,u(\cdot, t^*) \in L^2(|x|^8dx),\;\;\text{and}\;\;u(\cdot, t)\in L^2(|x|^{7-})- L^2(|x|^7),\;\;\forall 
\,t\notin\{ 0, t^*\},
$$
with $\,t^*$ as in \eqref{thetime}.
\vskip.1in
(b) The value of $t^*$ described in \eqref{thetime} can be motivated from the identity
$$
\frac{d\,\,}{dt}\,\int_{-\infty}^{\infty}x\,u(x,t)dx=\frac{1}{2} \|u(\cdot,t)\|_2^2=\frac{1}{2} \|u_0\|_2^2,
$$
(using the second conservation law which tells us that the $L^2$ norm of the real solution is preserved by the solution flow)
 which describes the time evolution of the first momentum of the solution
$$
\int_{-\infty}^{\infty}x\,u(x,t)dx= \int_{-\infty}^{\infty}x\,u_0(x)dx +\frac{t}{2} \|u_0\|_2^2.
$$
So assuming that
\begin{equation}
\label{condition}
\int_{-\infty}^{\infty}x\,u_0(x)dx\neq 0,
\end{equation}
one looks for the times where the average of the first momentum of the solution vanishes, i.e. for $t$ such that
$$
\int_0^t \int_{-\infty}^{\infty}x\,u(x,t')dx \,dt'= \int_0^t (\int_{-\infty}^{\infty}x\,u_0(x)dx +\frac{t'}{2} \|u_0\|_2^2) dt'=0,
$$
 which under the assumption \eqref{condition} has a unique solution $t=t^*$ given by the formula in \eqref{thetime}.
\vskip.1in
(c) To prove Theorem \ref{theorem7} we shall work with the integral equation version of the problem 
\eqref{BO}. Roughly speaking, from the result in \cite{MoSaTz}  one cannot regard the nonlinear term as a 
perturbation of the linear one. So to obtain our result we use an argument similar to that introduced in \cite{Io2}. 
This is based on the special structure of the equation and allows us 
to reduce the contribution of two terms in the integral 
equation to just one. Also the use of the integral equation in the proof and the result in \cite{MoSaTz} 
explains our assumption $u_0\in \dot Z_{5,4}$ instead of the expected one from the differential equation point of view 
$u_0\in \dot Z_{4,4}$.  
\vskip.1in
(d) One may ask if it is possible to have a stronger decay at $t=t^*$ than the one described in \eqref{attime}.
In this regard, our argument shows that for $u_0\in \dot Z_{6,5}$ it follows that
$$
u(\cdot,t^*)\in \dot Z_{5,5}
$$
if and only if
\begin{equation}\label{decaymore}
\int_0^{t^*} \,\int_{-\infty}^{\infty} x^2\,u(x,t)dx\,dt=0.
\end{equation}
However, the time evolution of the second momentum of the solution does not seem to have a simple expression which
allows to verify the identity \eqref{decaymore}.
 \vskip.1in
(e) The result in Theorem \ref{theorem7} can be extended to higher powers of the BO equations
$$
\partial_t u + \mathcal  H\partial_x^2u +u^{2k+1}\partial_x u = 0,\;\;\;\;k=0,1,2,...
$$
where the formula \eqref{thetime} for $t^*$ in this case will be given as the solution of the equation
$$
\int_0^{t^*}(\,\int x\,u_0(x)dx +\frac{1}{2k+2}\,\int_0^t\|u(t')\|_{2k+2}^{2k+2}\,dt')\,dt=0,\;\;\;\;k=1,2,...
$$
It is clear that if such a time $t^*$ exists it is unique.
\vskip.1in
(f) A close inspection of the proof of Theorem C in \cite{GFGP} gives us the following result which allows us to 
establish a 
uniqueness result with a condition involving only two times $t_1=0,\;t_2\neq 0$ for a suitable class of solutions:

 \begin{theorem}  \label{theorem7a} 
 Let $u\in C(\R : \dot Z_{7/2,3})$ be a solution of the IVP \eqref{BO}
 for which there exist two times $t_1,\,t_2\in \R,\,\,t_1\neq t_2\, $ such that
$$
u(\cdot,t_j)\in \dot Z_{7/2,7/2}.
$$
  If   
$$
\int x\,u(x,t_1)dx=0\;\;\;\;\;\;\;\text{or}\;\;\;\;\;\;\;\int x\,u(x,t_2)dx=0,
$$
then
$$ 
 u(x,t)\equiv 0.
$$
 
\end{theorem}
\vskip.1in
(g) In a forthcoming work we shall consider the extensions of the results established here and those
 in \cite{GFGP} to solutions of the IVP for the dispersive model
$$
\partial_t u + D_x^{1+a}\partial_xu +u\partial_x u = 0, \qquad t, x\,\in \R,\,\,\;a\in(0,1),
$$
where
$$
 D_x=(-\partial_x^2)^{1/2}=\mathcal H \partial_x.
$$
Thus, the cases $a=0$ and $a=1$ correspond to the BO equation and the KdV equation, respectively. 
\vskip.1in

 We recall that if  for a solution $u\in C(\R:H^s(\R)), \,s\geq 0$ of \eqref{BO} 
one has that $\,\exists \,t_0\in\R$ such that $u(x,t_0)\in H^{s'}(\R),\,s'>s$, 
then  $u\in C(\R:H^{s'}(\R))$. So the propagation of the $H^s(\R)$ regularity of the solution is not an issue.

\vskip.1in

As it was mentioned above, the proof of Theorem \ref{theorem7a} is contained in the proof of Theorem C given in 
\cite{GFGP},
therefore it will be omitted.
\vskip.1in

 The rest of this paper is organized as follows: section 2 contains all the estimates
 needed in the proof of Theorems \ref{theorem6}
 and \ref{theorem7}. The proof of Theorem \ref{theorem6} will be given in section 3.
 Theorem \ref{theorem7} will be proven in section 4. 
 
 \end{section}
 
\begin{section}{Preliminary Estimates}\label{S: 2}

 As in \cite{GFGP} we shall use the following generalization of the Calder\'on commutator estimates \cite{Ca} found in \cite{DaMcPo}:
 
   \begin{theorem}\label{theorem8}  
  For any $ \,p\in(1,\infty)$  and $l,\,m\in\mathbb  Z^+\cup\{0\},\,l+m\geq 1$  exists 
$\, c=c(p;l;m)>0$ such that
 \begin{equation}
\label{777}
 \| \partial_x^l[\mathcal H;\,a]\partial_x^m f\|_p\leq c \|\partial_x^{l+m} a\|_{\infty} \|f\|_p.
\end{equation}
\end{theorem}
\vskip.1in

We shall also use the pointwise identities
$$
[\mathcal H; x]\partial_x f = [\mathcal H; x^2]\partial_x^2 f=0,
$$
and more generally
$$
[\mathcal H; x] f =0\;\;\;\;\text{if and only if }\;\;\;\;\int fdx=0,
$$
 and
$$
[\mathcal H; x^2] f =0\;\;\;\;\text{if and only if }\;\;\;\;\int fdx=\int \,xf dx=0.
$$ 
\vskip.1in
   
  To justify the finiteness of the quantities involved in the energy estimate used in the proof of Theorem \ref{theorem6} we introduce the  truncated weights $w_N(x)$.  
  Using the notation  $\langle x\rangle =(1+x^2)^{1/2}$ we define
\begin{equation}
\label{truncw}
 w_N(x)=
\begin{cases}
\langle x
\rangle& \text{if\,\,$|x|\le N$,}\\
2N  &  \text{if\,\,$|x|\ge 3N,$}
\end{cases}
\end{equation}
 $w_N(x) $ smooth and non-decreasing in $|x|$ with  $w_{N}'(x)\leq 1$ for all $x\geq 0$.
 We observe that
  $$
  x\,w'_N(x)\leq c w_N(x),
  $$
  where the constant $c$ is independent of $N$.

\end{section}
  
 \begin{section}{Proof of Theorem \ref{theorem6}}\label{S: 3}
 
We take  two solutions $u_1,\,u_2$ of \eqref{BO} whose  data $\,u_{1,0},\;u_{2,0}\,$ satisfy
  \begin{equation}
  \label{dato1}
  u_1(x,0)=u_{1,0}(x),\,\;u_2(x,0)=u_{2,0}(x)\in Z_{4,4},
  \end{equation}
  with
  \begin{equation}
  \label{dato2}
  \begin{cases}
  \begin{aligned}
   \;\;\;\;\; \int_{-\infty}^{\infty} u_1(x,0)\,dx&=\int_{-\infty}^{\infty} u_2(x,0)\,dx,\\
    \\
    \;\;\;\; \int_{-\infty}^{\infty}\,x\, u_1(x,0)\,dx&=\int_{-\infty}^{\infty} \,x\,u_2(x,0)\,dx,\\
    \\
 \;\;\;\; \|u_{1,0}\|_2=\|u_{2,0}\|_2&,\;\;\;\;\;\;\;\;\;u_{1,0}\neq u_{2,0}\\
\\
 \;\;\;\;u_{1,0}\neq 0\,,&\,\,\,\,\,\,\,\,\,\,\,\,u_{2,0}\neq 0.
  \end{aligned}
  \end{cases}
  \end{equation}
  \vskip.1in
   Thus, from the result in \cite{GFGP} it follows that
  $$
  u_1,\;u_2\in C(\R: Z_{4,5/2^-}).
  $$
 Defining
 $$
 v(x,t)=u_1(x,t)-u_2(x,t),
 $$
 one sees that $v$ verifies the linear equation
 \begin{equation}
 \label{eqv}
 \partial_t v + \mathcal  H\partial_x^2v +u_1\partial_x v +  \partial_xu_2 v= 0,
 \end{equation}
with
 \begin{equation}
 \label{diff}
 v\in C(\R: Z_{4,5/2^-}),
 \end{equation}
 and 
 \begin{equation}
 \label{mom}
  \int_{-\infty}^{\infty} v(x,t)\,dx= \int_{-\infty}^{\infty} x \,v(x,t)\,dx=0,\;\;\;\;\forall t\in \R.
 \end{equation}
The identities in \eqref{mom}  follow by combining  our hypothesis \eqref{dato2}, the first conservation law 
$$
 \int_{-\infty}^{\infty}\,u_j(x,t)\,dx =  \int_{-\infty}^{\infty}\,u_{j,0}(x)\,dx ,\;\;\;\;\;\;\;\forall \,t\in\R,\;\;j=1,2,
 $$
and the identity
  $$
  \frac{d\;\;}{dt}  \int_{-\infty}^{\infty} x\,u_j(x,t)\,dx =\frac{1}{2} \|u_j(t)\|_2^2=\frac{1}{2}\|u_{j,0}\|_2^2,\;\;\;\;\;\;\;\forall \,t\in\R,\;\;\;j=1,2.
  $$
  
Now, differentiating the equation in \eqref{eqv} and multiplying the result by $w_N^2$ we get
  \begin{equation}
 \label{eqv1}
 \partial_t(w_N^2\partial_x v) +  w_N^2\mathcal  H\partial_x^2\,\partial_x v + w_N^2\,\partial_x(u_1\partial_x v + v \partial_xu_2)= 0.
 \end{equation}
We rewrite the second term in \eqref{eqv1} as
\begin{equation}
\label{a1}
\begin{aligned}
&w_N^2\mathcal  H\partial_x^2\,\partial_x v\\
&=\mathcal H (w_N^2 \partial_x^2\,\partial_x v) - [\mathcal H; w_N^2]\partial_x^3\,v\\
&=\mathcal H \partial_x^2(w_N^2 \,\partial_x v) -2 \mathcal H(\partial_xw_N^2 \partial_x^2v)-\mathcal H(\partial_x^2w_N^2 \,\partial_xv)- [\mathcal H; w_N^2]\partial_x^3\,v\\
&=G_1+G_2+G_3+G_4.
\end{aligned}
\end{equation}
 
 Theorem \ref{theorem8} yields the inequality
 $$
 \|G_4\|_2=\| [\mathcal H; w_N^2]\partial_x^3\,v\|_2  \leq c\,\|v\|_2,
 $$
 with $c$ denoting a constant independent of $N$ which may change from line to line. Also one has that
 $$
 \|G_3\|_2=\|\mathcal H(\partial_x^2w_N^2 \,\partial_xv)\|_2\leq c \,\|\partial_x v\|_2.
 $$
 To control $\|G_2\|_2$ we use integration by parts to get that
 $$
 \|w_N\,\partial_x^2v\|_2^2\leq \|w_N^2\,\partial_xv\|\,\|\partial_x^3v\|_2+\|\partial_xv\|_2^2,
 $$
 so
 $$
 \|G_2\|_2=\|\mathcal H(\partial_xw_N^2 \,\partial^2_xv)\|_2\leq \|w_N  \,\partial^2_xv\|_2
\leq c(\|w_N^2\partial_x v\|_2+\|v(t)\|_{3,2}).
 $$
In the energy estimate the contribution  of the term $G_1$ is null, since inserting it  in \eqref{eqv1},
  multiplying the equation \eqref{eqv} by $w_N^2\partial_x v$, and integrating
 the result in the space variable after integration by parts it vanishes. It remains 
to bound the contribution from the third term in \eqref{eqv1} in the energy estimate, i.e.
 $$
N_1(t)=| \int_{-\infty}^{\infty}w_N^2\,\partial_x(u_1\partial_x v + v \partial_xu_2)\,w_N^2\partial_x v\,dx |.
$$
Using the hypotheses and integration by parts it follows that for any $T>0$
$$
\aligned
N_1(t)&\leq c_T (\|\partial_x u_1(t)\|_{\infty}+\|\partial_x u_2(t)\|_{\infty}) \|w_N^2\partial_x v\|_2^2 \\
&\;\;\;\;\;+ c_T\|u_1\|_{\infty}\|x\partial_xv\|_2\|w_N^2\partial_x v\|_2+ c_T \|\partial_x^2u_2(t)\|_{\infty} \| x^2v\|_2 \|w_N^2\partial_x v\|_2\\
&\leq c_T(\|w_N^2\partial_x v(t)\|_2+\|w_N^2\partial_x v(t)\|_2^2),\;\;\;\;\;\;\;\forall t\in [-T,T],
\endaligned
$$
with $c_T$ denoting a constant  depending on the initial solutions and on their data but independent of 
$N$ whose value may change from to line.  
Collecting the above information we conclude that for any $T>0$ 
$$
\sup_{t\in[-T,T]} \| w_N^2\,\partial_x v(t)\|_2 < c_T.
$$
Therefore, taking $N\uparrow \infty$ it follows that for any $T>0$
\begin{equation}
\label{step1}
\sup_{t\in[-T,T]} \| x^2\,\partial_x v(t)\|_2 < M_T,
\end{equation}
with $M_T$ denoting a constant depending only on initial parameters and on $T$ and whose value may change from  line to line.
Since by hypothesis we have
$$
\sup_{t\in[-T,T]} \| \partial^3_x v(t)\|_2 <  M_T,
$$
 by integration by parts one gets that for any $T>0$ 
 $$
 \sup_{t\in[-T,T]} \|x \,\partial^2_x v(t)\|_2 < M_T.
$$

Next, from  the identity
 $$
 x\,\mathcal H\,\partial_x^2 v =\mathcal H\,(x \,\partial_x^2v) =\mathcal H \partial_x^2(xv)-2\mathcal H \partial_x v,
 $$
we get the equation for $w_N^2 x v$
\begin{equation}
 \label{eqv2}
 \begin{aligned}
& \partial_t (w_N^2 x v) + \mathcal  H\partial_x^2(w_N^2 x v) - 2\mathcal H(\partial_x w_N^2\partial_x(xv))-\mathcal H(\partial_x^2w_N^2\,xv)\\
 & -[\mathcal H; w_N^2]\partial_x^2(xv) -2 w_N^2 \mathcal H\partial_x v+
 w_N^2\,x(u_1\partial_x v + v \partial_xu_2)= 0.
\end{aligned}
 \end{equation}
  
We recall that  for all $\,t\in\R$
   \begin{equation}
    \label{ke} 
\int_{-\infty}^{\infty}v(x,t)\,dx=\int_{-\infty}^{\infty}x\,v(x,t)\,dx=0,
   \end{equation}
 so that
 \begin{equation}
 \label{key}
 xH(v)=H(xv),\;\;\;\;\text{and}\;\;\;\;x^2H(v)=H(x^2v).  
 \end{equation}
  
   The following string of estimates
  $$
  \| \mathcal H(\partial_x w_N^2\partial_x(xv))\|_2\leq c(\|w_N x \partial_xv\|_2 +\| w_N v\|_2)\leq c(\| x^2\partial_x v\|_2+\| xv\|_2),
  $$
  $$
 \|  \mathcal H(\partial_x^2w_N^2\,xv)\|_2\leq c\| xv\|_2,
 $$
 (by Theorem \ref{theorem8})
 $$
 \| [\mathcal H; w_N^2]\partial_x^2(xv)\|_2\leq c \|\partial_x^2w_N^2\|_{\infty}\|x v\|_2\leq c\|x v\|_2,
 $$
 (by \eqref{key})
 $$
 \aligned
& \| w_N^2 \mathcal H\partial_x v\|_2\leq \| (1+x^2) \mathcal H\partial_x v\|_2\\
&\leq \| \mathcal \partial_x v\|_2 + \| x^2\mathcal H\partial_x v\|_2
\leq \|\partial_x v\|_2+ \| x \mathcal H(x\partial_x v)\|_2\\
&\leq \|\partial_x v\|_2+ \| \mathcal H(x^2\partial_x v)\|_2\leq  \|\partial_x v\|_2+ \| x^2\partial_x v\|_2,
\endaligned
$$
and (by integrating by parts)
 $$
 \aligned
& |\int w_N^2\,x(u_1\partial_x v + v \partial_xu_2) w_N^2 x v dx|\\
&\leq 
 (\|\partial_x u_1\|_{\infty} + \|\partial_x u_2\|_{\infty}) \| w_N^2 xv\|_2^2 +\| u_1\|_{\infty}\|x^2 v\|_2  \| w_N^2 xv\|_2.
 \endaligned
 $$
  inserted in  the energy estimate for \eqref{eqv2}  
  together with the  result in the previous step \eqref{step1} allows us to conclude that for any $T>0$
  $$
  \sup_{t\in[-T,T]} \|w_N^2 xv\|_2\leq c_T,
$$
  with $c_T$ independent of $N$. Hence, it follows that
 \begin{equation}
  \label{step3}
  \sup_{t\in[-T,T]} \| x^3\,v\|_2\leq M_T.
  \end{equation} 

Now, we shall estimate  $w_N^2x\partial_x v$. From the equation
$$
\partial_t(x\partial_x v) +  \mathcal  H\partial_x^2\,(x\partial_x v) - 2\mathcal H\partial_x^2 v + x\,\partial_x(u_1\partial_x v + v \partial_xu_2)= 0,
$$
we obtain that 
  \begin{equation}
 \label{eqv3}
 \begin{aligned}
&\partial_t(w_N^2 x\partial_x v) +  \mathcal  H\partial_x^2\,(w_N^2 x\partial_x v) 
-2\mathcal H(\partial_xw_N^2  \partial_x(x\partial_xv)) \\
&- \mathcal H(\partial^2_x(w_N^2) x\partial_xv)
-[\mathcal H;w_N^2]\partial_x^2(x\partial_xv) -2w_N^2\mathcal H\partial_x^2 v \\
&+ w_N^2x\,\partial_x(u_1\partial_x v + v \partial_xu_2)= 0.
\end{aligned}
\end{equation}

 By integration by parts one gets that
 \begin{equation}
 \label{int1}
 \begin{aligned}
& \|w_N\partial_x^3v\|_2^2\leq c(\|w_N^2\partial_x^2v\|_2\,\|\partial_x^4v\|_2 + \|\partial_x^2v\|_2^2),\\
& \|w_N^2\partial_x^2v\|_2^2\leq c(\|w_N^3\partial_xv\|_2\,\|w_N\partial_x^3v\|_2 + \|w_N\partial_xv\|_2^2),
\end{aligned}
\end{equation}
with a constant $c$ independent of $N$.
We observe  that for each $N$ fixed all the quantities in \eqref{int1} are finite. Hence, from \eqref{int1} it follows that
 \begin{equation}
 \label{int2}
 \begin{aligned}
& \|w_N\partial_x^3v\|_2\leq c(\|w_N^3\partial_x v\|^{1/3}_2\,\|\partial_x^4v\|^{2/3}_2 + \|w_N\partial_xv\|_2+\|v\|_{4,2}),\\
& \|w_N^2\partial_x^2v\|_2\leq c(\|w_N^3\partial_xv\|^{2/3}_2\,\|\partial_x^4v\|^{1/3}_2 + \|w_N^3\partial_xv\|_2 +\|v\|_{2,2}),
\end{aligned}
\end{equation}
with $c$ independent of $N$.

Returning to the equation \eqref{eqv3} we shall use Theorem \ref{theorem8} to get that
 $$
 \aligned
& \|w_N^2\mathcal H\partial_x^2 v\|_2\\
&\leq \|\mathcal H(w_N^2\partial_x^2 v)\|_2 +\| [\mathcal H;w_N^2]\partial_x^2v\|_2\\
&\leq c(\|w_N^2\partial_x^2 v\|_2 + \|\partial_x^2w_N^2\|_{\infty}\|v\|_2)=c(D_1+\|v\|_2).
\endaligned
$$
Thus, by combining the second inequality in \eqref{int2} and Young's inequality it follows that
\begin{equation}
 \label{D1}
D_1=\|w_N^2\partial_x^2 v\|_2\leq c(\|w_N^2 x \partial_xv\|_2+\|v\|_{4,2} ).
\end{equation}
 Theorem \ref{theorem8} yields the inequality
 $$
 \| [\mathcal H; w_N^2]\partial_x^2(x\partial_xv)\|_2  \leq c\,\|x\partial_xv\|_2,
 $$
 with $c$ independent of $N$ whose value may change from line to line. Also one has
 $$
\|\mathcal H(\partial_x^2w_N^2 \,x\partial_xv)\|_2\leq c \,\|x \partial_x v\|_2.
 $$
 To control the third term in \eqref{eqv3} we write
 $$
 \aligned
& \|\partial_x(w_N^2)\,\partial_x(x\partial_xv)\|_2\\
&\leq c(\|w_N w_N'x\partial_x^2v\|_2+\|w_Nw_N'\partial_xv\|_2\\
&\leq c(\|w_N^2\partial_x^2v\|_2 +\|x\partial_x v\|_2)=c(D_1+\|x\partial_x v\|_2),
\endaligned
 $$
 with $D_1$ as in \eqref{D1}.
So besides the first two terms in \eqref{eqv3} it remains to bound the contribution from the last term 
in the energy estimate, i.e.
 $$
N_2(t)=| \int_{-\infty}^{\infty}w_N^2\,x\partial_x(u_1\partial_x v + v \partial_xu_2)w_N^2\,x\partial_x v\,dx |.
$$
Using our hypothesis and integration by parts it follows that for any $T>0$
$$
\aligned
N_2(t)&\leq c(\|\partial_xu_1\|_{\infty} +\|\partial_xu_2\|_{\infty})\|\,w_N^2x\partial_xv\|_2^2 \\
&+(\|u_1\|_{\infty}\|x^2 \partial_xv\|_2+\|\partial_x^2u_2\|_{\infty}\|x^3v\|_2)\|\,w_N^2x\partial_xv\|_2,
\endaligned
$$
with $c_T$ depending on the initial solutions and their data but independent of $N$.  Collecting 
the above information we conclude that for any $T>0$ 
$$
\sup_{t\in[-T,T]} \| w_N^2\,x\partial_x v(t)\|_2 < c_T,
$$
with $c_T$ depending  on the initial solutions $u_1,\,u_2$, the initial data, and on $T$, but independent of $N$. Therefore, taking $N\uparrow \infty$ it follows that for any $T>0$
\begin{equation}
\label{step11}
\sup_{t\in[-T,T]} \| x^3\,\partial_x v(t)\|_2 < M_T,
\end{equation}
with $M_T$ denoting a generic constant which may change line to line but depending only on initial parameters and on $T$.
From \eqref{int2} we have
$$
\sup_{t\in[-T,T]} \| x^2\partial^2_x v(t)\|_2 <  M_T,
$$
 by integration by parts one gets that for any $T>0$ 
 \begin{equation}
 \label{previous2}
 \sup_{t\in[-T,T]} \|x \,\partial^3_x v(t)\|_2 < M_T.
\end{equation}

Using  the identity
 $$
 x^2\,\mathcal H\,\partial_x^2 v =\mathcal H \partial_x^2(x^2v)-4\mathcal H \partial_x(x v) + 2\mathcal Hv,
 $$
we get the equation for $w_N^2 x^2 v$
\begin{equation}
 \label{eqv22}
 \begin{aligned}
& \partial_t (w_N^2 x^2 v) + \mathcal  H \partial_x^2(w_N^2 x^2 v) - 2\mathcal H(\partial_x w_N^2\partial_x(x^2v))\\
&-\mathcal H(\partial_x^2w_N^2\,x^2v)-[\mathcal H; w_N^2]\partial_x^2(x^2v) -4 w_N^2 \mathcal H \partial_x(xv)\\
&+2w_N^2\mathcal H v+w_N^2\,x^2(u_1\partial_x v + v \partial_xu_2)= 0.
\end{aligned}
 \end{equation}
  
We shall use \eqref{ke} and \eqref{key} and deduce the following  estimates: \newline
   (using \eqref{key})
  $$
  \| w_N^2 \mathcal Hv \|_2\leq \|(1+x^2)\mathcal Hv\|_2 \leq \|xv\|_2 +\| x^2 Hv\|_2\leq \|xv\|_2 +\| x^2 v\|_2,
  $$
  (using \eqref{key} and \eqref{previous2})
  $$
  \aligned
& \| w_N^2 \mathcal H\partial_x(xv)\|_2\leq \| (1+x^2)\mathcal H\partial_x(xv)  \|_2\\
&\leq  \| \mathcal H\partial_x(xv)  \|_2 +  \| x^2\mathcal H\partial_x(xv)  \|_2\\
& \leq \| \partial_x(xv)  \|_2 +  \| x\mathcal H(x\partial_x(xv))  \|_2\\
&\leq   \| \partial_x(xv)  \|_2 +  \| x\mathcal H(\partial_x(xv))  \|_2 + 2 \| x\mathcal H(xv)  \|_2\\
&\leq   \| \partial_x(xv)  \|_2 +  \| x \partial_x(xv)  \|_2 + 2 \| x^2v  \|_2,
\endaligned
 $$
 (using Theorem \ref{theorem8})
 $$
 \| [\mathcal H; w_N^2]\partial_x^2(x^2v)\|_2\leq c \|\partial_x^2w_N^2\|_{\infty}\|x^2 v\|_2\leq c\|x^2 v\|_2,
 $$
 $$
 \| \mathcal H (\partial_x^2w_N^2 x^2 v)\|_2\leq \|x^2 v\|_2,
$$
(using \eqref{step11})
$$
\aligned
& \| \mathcal H (\partial_xw_N^2 \partial_x(x^2 v))\|_2\leq \| \partial_xw_N^2 \partial_x(x^2 v)  \|_2
\\
& \leq 8(\|w_Nw_N'xv\|_2+\|w_Nw_N' x^2\partial_xv\|)
\leq 8(\| x^2v\|_2 +\|w_N^3\partial_xv\|_2),
\endaligned
$$
and integrating by parts (for the last term in \eqref{eqv22})
 $$
 \aligned
& |\int w_N^2\,x^2(u_1\partial_x v + v \partial_xu_2) w_N^2 x^2 v dx|\\
&\leq 
 (\|\partial_x u_1\|_{\infty} + \|\partial_x u_2\|_{\infty}) \| w_N^2 x^2v\|_2^2 +\| u_1\|_{\infty}\|w_N^2 x v\|_2  \| w_N^2 x^2v\|_2.
 \endaligned
 $$
 Collecting this information  in  the energy estimate for \eqref{eqv22}  
  together with the  result in the previous steps \eqref{previous2} and  \eqref{step11} allows us to conclude that for any $T>0$
  $$
  \sup_{t\in[-T,T]} \|w_N^2 x^2v\|_2\leq M_T,
$$
  with $M_T$ independent of $N$. Hence, it follows that
 \begin{equation}
  \label{step133}
  \sup_{t\in[-T,T]} \| x^4\,v\|_2\leq M_T.
  \end{equation} 
 Hence, for any $T>0$
 $$
v\in L^{\infty}([-T,T]: Z_{4,4}),
$$
which yields the desired result.

 \end{section}


\begin{section}{Proof of Theorem \ref{theorem7}}\label{S: 4}

We introduce the notation
\begin{equation}
\label{defF}
F_j(t,\xi, \widehat{u}_0)=\partial_{\xi}^j(e^{-it|\xi|\xi}\widehat{u}_0),\;\;\;\;j=0,1,2,3,4.
\end{equation}

Therefore, 
\begin{equation}
 \begin{split}
F_3(t,\xi,\widehat{u}_0)&=\partial_{\xi}^3(e^{-it|\xi|\xi}\widehat{u}_0)\\
&=e^{-it|\xi|\xi}(8it^3\xi^3\widehat{u}_0-12t^{2}\xi\widehat{u}_0-12t^{2}\xi^2\partial_{\xi} \widehat{u}_0\\
&\,\,\,\,-6it\,\text{sgn}(\xi)\partial_{\xi} \widehat{u}_0-6it|\xi|\partial_{\xi} ^2\widehat{u}_0
-2it\delta\widehat{u}_0+\partial_{\xi} ^3\widehat{u}_0).
\end{split}
\label{F3}
\end{equation}

 We observe that  since the initial data $u_0$ has zero mean value (i.e. $\,\widehat u_0(0)=0$) the term involving the Dirac delta in  
\eqref{F3} vanishes.  
Thus, under the assumption that $u_0$ has zero mean value one finds that
\begin{equation}
 \begin{split}
F_4(t,\xi,\widehat{u}_0)&=\partial_{\xi}^4(e^{-it|\xi|\xi}\widehat{u}_0)\\
&=e^{-it|\xi|\xi}(12 t^2 \widehat{u}_0+48 i t^3\xi |\xi| \widehat{u}_0+ 16 t^4 \xi^4 \widehat{u}_0\\
&\,\,\,\,\,-48 t^2 \xi\partial_{\xi} \widehat{u}_0-6it  \delta\partial_{\xi} \widehat{u}_0 +24 i t^3 |\xi|\xi^2\partial_{\xi} \widehat{u}_0 \\
&\,\,\,\,-12it \,\text{sgn}(\xi)\partial_{\xi}^2\widehat{u}_0-24 t^2\xi^2\partial_{\xi}^2\widehat{u}_0-8it|\xi|\partial_{\xi}^3\widehat{u}_0+\partial_{\xi} ^4\widehat{u}_0)\\
&=E_1(t,\xi,\widehat{u}_0)+.....+E_{10}(t,\xi,\widehat{u}_0).
\end{split}
\label{F4}
\end{equation}

Hence 
 \begin{equation}
\widehat{u}(\xi,t)= F_0(t,\xi,\widehat{u}_0)-\int_0 ^tF_0(t-t',\xi,\widehat{z}(t'))\,dt',
\label{ec}
\end{equation}
and
\begin{equation}
\partial_{\xi}^4\widehat{u}(\xi,t)= F_4(t,\xi,\widehat{u}_0)-\int_0 ^tF_4(t-t',\xi,\widehat{z}(t'))\,dt',
\label{ec4}
\end{equation}
where 
$$
\widehat{z}= \frac{1}{2}\widehat {\partial_x u^2}=i \frac{\xi}{2} \widehat{u}\ast\widehat{u}.
$$
\vskip.1in
We shall see  that if $u_0\in  \dot Z_{4,4}$ all terms $\,E_j,\,j=1,...,10$ in \eqref{F4} except 
\begin{equation}
\label{E5}
E_5(t,\xi,\widehat{u}_0)=e^{-it|\xi|\xi}(-6it  \delta\partial_{\xi} \widehat{u}_0(\xi))=-6it  \delta\partial_{\xi} \widehat{u}_0(0)=-6\,t\,\delta\int x\,u_0(x)dx,
\end{equation}
are in $L^2(\R)$. Thus, we have
\begin{equation}
\label{L2}
\begin{aligned}
\begin{cases}
\;\|E_1\|_2&=\|12 t^2\,  e^{-it|\xi|\xi}\widehat u_0\|_2 \leq c_t\,\|u_0\|_2,\\
\;\|E_2\|_2&=\| 48 t^3\,|\xi|\xi e^{-it|\xi|\xi}\widehat u_0\|_2 \leq c_t\, \|\partial^2_xu_0\|_2,\\
\;\|E_3\|_2&=\| 16 t^4\,\xi^4e^{-it|\xi|\xi}\widehat u_0\|_2 \leq c_t\, \|\partial^4_xu_0\|_2,\\
\;\|E_4\|_2&=\| 48 t^2\,\xi e^{-it|\xi|\xi}\partial_{\xi}\widehat u_0\|_2 \leq c_t\, (\|u_0\|_2+\|x\partial_xu_0\|_2),  \\
\;\|E_6\|_2&=\| 24 t^3\,|\xi|\xi^2e^{-it|\xi|\xi}\partial_{\xi}\widehat u_0\|_2 \leq c_t\,(\|x\partial_x^3u_0\|_2+\|\partial_x^2u_0\|_2),\\
\;\|E_7\|_2&=\| 24 t^2\,\xi^2e^{-it|\xi|\xi}\partial_{\xi}^2\widehat u_0\|_2 \leq c_t\,(\|u_0\|_2+\|x^2\partial_x^2u_0\|_2+\|x\partial_xu_0\|_2),\\
\;\|E_8\|_2&=\| 12 t\,\text{sgn}(\xi)e^{-it|\xi|\xi}\partial_{\xi}^2\widehat u_0\|_2 \leq c_t\,\|x^2u_0\|_2,\\
\;\|E_9\|_2&= \| 8 t \,|\xi|e^{-it|\xi|\xi}\partial_{\xi}^3\widehat u_0\|_2 \leq c_t\,(\|x^3\partial_xu_0\|_2 +\|x^2u_0\|_2),\\
\;\|E_{10}\|_2&=\|e^{-it|\xi|\xi}\partial_{\xi}^4\widehat u_0\|_2 \leq c_t\,\|x^4u_0\|_2.
\end{cases}
\end{aligned}
\end{equation}

Since $u_0\in Z_{4,4}=H^4(\R)\cap L^2(|x|^8dx)$, by using interpolation it follows directly that all the terms 
on the left hand side 
of \eqref{L2} are bounded. In fact,
$$
E_j\in C(\R:L^2(\R)),\;\;\;\text{for}\,\,\,\;1\leq j\leq 10,\,\,j\neq 5.
$$
Thus,
\begin{equation}
\label{abc} 
F_4-E_5= F_4(t,\xi,\widehat{u}_0)-E_5(t,\xi,\widehat{u}_0)\in C(\R:L^2(\R)),
\end{equation}
and since $\,E_5\in C(\R:H^{-(1/2+\epsilon)}(\R)),\,\epsilon>0,$ one has
$$
F_4\in C(\R:H^{-(1/2+\epsilon)}(\R)),\;\;\;\epsilon>0.
$$

 Next, we shall consider the integral term in \eqref{ec4}
 $$
\Omega=\Omega(t,\xi,\widehat{z})\equiv \int_0 ^tF_4(t-t',\xi,\widehat{z}(t'))\,dt' ,
$$
with 
$$
\widehat{z}= \frac{1}{2}\widehat {\partial_x u^2}=i \frac{\xi}{2} \widehat{u}\ast\widehat{u}.
$$

We are assuming that $\,u_0\in \dot Z_{5,4}$, therefore from Theorem A we have that the 
corresponding solution $u(x,t)$ of \eqref{BO} satisfies that
$$
u\in C(\R:\dot Z_{5,7/2-}).
$$
Thus, 
$$
\widehat u\in C(\R: Z_{7/2-, 5}),
$$
and hence
$$
\widehat {u}\ast\widehat {u}\in C(\R: Z_{6, 5}),
$$
so we can conclude that
\begin{equation}
\label{aaa}
\xi\,\widehat{u}\ast\widehat{u} \in C(\R: Z_{4,4}).
\end{equation}

Above we have seen that if $u_0\in \dot Z_{4.4}$, then  all the ten terms of  
$F_4(t,\xi,\widehat u_0)$ (see \eqref{F4}) except $E_5$  (see \eqref{E5}) are in $C(\R:L^2(\R))$.
The same argument, the fact that $u\partial_xu$ has mean value zero, \eqref{L2}, and \eqref{aaa}  proves
that 
$$
\aligned
\Omega(t,\xi,\widehat{z})&= \int_0 ^tF_4(t-t',\xi,\widehat{z}(t'))dt'\\
\\
&=\sum_{j=1}^{10}\int_0^t\,E_j(t-t',\xi,\widehat{z}(t'))dt'\equiv \sum_{j=1}^{10} \,B_j(t,\xi,\widehat{z}),
\endaligned
$$
with
$$
B_j=B_j(t,\xi,\widehat{z})\in C(\R:L^2(\R)),\,\,\,1\leq j\leq 10,\,\,j\neq 5.
$$
Therefore, 
\begin{equation}
\label{nonL2}
\begin{aligned}
& \Omega(t,\xi,\widehat{z})-\int_0^t\,E_5(t-t',\xi,\widehat{z}(t'))dt'\\
&= \Omega+ 6i\int_0^t(t-t')e^{-i(t-t')|\xi|\xi}\delta \partial_{\xi}
(\frac{i \xi}{2} \widehat{u}\ast\widehat{u})(\xi,t') dt' \\
&=\Omega + 6i\,\delta \int_0^t(t-t')\, \partial_{\xi}(\frac{i \xi}{2} \widehat{u}\ast\widehat{u})(0,t') dt'\\
&\equiv  \Omega - B_5\in C(\R:L^2(\R)).
 \end{aligned}
\end{equation}
We observe that
\begin{equation}
\label{001}
\begin{aligned}
&\partial_{\xi}(\frac{i \xi}{2} \widehat{u}\ast\widehat{u})(0,t')=\widehat{-ixu\partial_xu}(0,t')=-i\;
\int_{-\infty}^{\infty}x\,u\partial_xu(x,t')dx\\
&=\frac{i}{2}\|u(t')\|_2^2=\frac{i}{2}\|u_0\|_2^2=i\,\frac{d\;\;}{dt}\int_{-\infty}^{\infty}x\, u(x,t)\,dx.
\end{aligned}
\end{equation}
 
Using \eqref{001} and integration by parts it follows that

\begin{equation}
\label{all}
\begin{aligned}
B_5&= - 6i\,\delta \int_0^t(t-t')\, \partial_{\xi}(\frac{i \xi}{2} \widehat{u}\ast\widehat{u})(0,t') dt'\\
&= - 6i\,\delta \int_0^t(t-t')\, (i\,\frac{d\;\;}{dt}\int x\,u(x,t')dx)\, dt'\\
&= 6\, \delta ( \;(t-t')\,\int x\, u(x,t')dx |_{t'=0}^{t'=t}\;+\;\int_0^t\,(\int x\, u(x,t')dx)\,dt')\\
&=- 6\,\delta(\,t\,\int x\,u_0(x)dx \,-\,\;\int_0^t\,(\int x\, u(x,t')dx)\,dt').
\end{aligned}
\end{equation}
Since, $\,B_5\in C(\R:H^{-(1/2+\epsilon)}(\R)),\,\;\epsilon>0$, one has that
$$
\Omega \in C(\R:H^{-(1/2+\epsilon)}(\R)),\,\epsilon>0.
$$
Collecting the above information we have that
$$
\partial_{\xi}^4\widehat{u}=F_4-\Omega\in C(\R:H^{-(1/2+\epsilon)}(\R)),\,\epsilon>0,
$$
and by \eqref{abc} and \eqref{nonL2}
$$
\aligned
&\partial_{\xi}^4\widehat{u}-E_5+B_5\\
&\;= F_4(t,\xi,\widehat{u}_0)-\int_0 ^tF_4(t-t',\xi,\widehat{z}(t'))\,dt' - E_5+B_5\\
&\;= F_4(t,\xi,\widehat{u}_0)-\Omega(t,\xi,\widehat{z}) - E_5+B_5\in C(\R:L^2(\R)).
\endaligned
$$
Now, using \eqref{001} it follows that
$$
\aligned
&E_5(t)-B_5(t)\\
&=-6\,t\,\delta\int x\,u_0(x)dx +
6\,\delta \Big(t\,\int x\,u_0(x)dx \,-\,\;\int_0^t\,(\int x \,u(x,t')dx\,)\,dt'\Big)\\
&=\,- 6\,\delta \,\int_0^t\,(\int x\, u(x,t')dx)\,dt')= - 6\,\delta\, \int_0^t ( \int x\,u_0(x)dx+\frac{t'}{2}\,\|u_0\|^2_2)dt'\\
&=\, - 6\,\delta \,(t\,\int x\, u_0(x)dx + \frac{t^2}{4}\|u_0\|_2^2),
\endaligned
$$
which vanishes only at
$$
t^*=\,-\, \frac{4}{\|u_0\|_2^2}\,\int_{-\infty}^{\infty} x\,u_0(x)dx,
$$
and at that time we have that
$$
\partial_{\xi}^4u(\cdot,t^*)\in L^2(\R),
$$
which yields the desired result.

\end{section}

\vskip.2in
{\underbar{ACKNOWLEDGMENT}}  :  
This work was done while G. P. was visiting the Department of Mathematics at Universidad Autonoma de 
Madrid-Spain whose hospitality he gratefully acknowledges. F.L. was partially supported by CNPq and FAPERJ.
G.P. was supported by NSF grant DMS-0800967.

\end{document}